\documentclass{article}

\usepackage{amsmath,amssymb,graphicx,color,subfig,hyperref}

\usepackage{amsfonts}
\usepackage{textcomp}
\usepackage{mathrsfs}
\usepackage{float}
\usepackage{epsfig}

\usepackage{enumitem}
\usepackage{verbatim}

\newcommand{\NN}{\mathbb{N}}

\newcommand{\RR}{\mathbb{R}}

\newtheorem{remark}{Remark}[section]

\newtheorem{definition}{Definition}[section]
\newtheorem{theorem}{Theorem}[section]
\begin{document}

	
	
	\newcommand{\Eqref}[1]{(\ref{#1})}

	

	\begin{center}
		\large{ \textbf{ {\Large Robust approximation algorithms for the \\
					\vskip 0.05cm
					detection of attraction basins in dynamical 
						systems}}}
	\end{center}

	\begin{center}
		Roberto Cavoretto$^{*}$, Alessandra De Rossi$^{*}$, Emma Perracchione$^{*}$ and \\   Ezio Venturino$^{*}$
	\end{center}
	
	\begin{center}
		$^{*}$Department of Mathematics G. Peano, University of Turin - Italy
	\end{center}
	\vskip 0.5cm
	
\textbf{Abstract.}	
In dynamical systems saddle points partition the domain into basins of attractions of the remaining locally stable equilibria. This problem is rather common especially in population dynamics models. Precisely, a particular solution of a dynamical system is completely determined by its initial condition and by the parameters involved in the model. Furthermore, when the omega limit set reduces to a point, the trajectory of the solution evolves towards the steady state. But, in case of multi-stability it is possible that several steady states originate from the same parameter set. Thus, in these cases the importance of accurately reconstruct the attraction basins follows. In this paper we focus on dynamical systems of ordinary differential equations presenting three stable equilibia and we design algorithms for the detection of the points lying on the manifolds determining the basins of attraction and for the reconstruction of such manifolds. The latter are reconstructed by means of the implicit partition of unity method which makes use of radial basis functions (RBFs) as local approximants. Extensive numerical test, carried out with a \textsc{Matlab} package  made available to the scientific community, support our findings.

\section{Introduction}
Mathematical modelling is applied in  major disciplines, such as biology, medicine and social sciences. The aim of such models lies in the prediction of the temporal evolution of the considered quantities (populations, cancer, divorces), \cite{Arrowsmith,Murray}. 

In a model involving a  set of ordinary differential equations,  a particular solution of the system is completely determined by the initial condition.  The latter establishes the steady state of the solution. Moreover, we may also note that, under some conditions imposed on the model parameters and depending on the initial state, the trajectories, i.e. the model solutions, at the end of the observation period may stabilize at different equilibria.
So the phase state of the dynamical system is partitioned into different regions, called the basins of attraction, depending on where the trajectories originating in them will ultimately stabilize. Thus, the prediction of a mathematical model  depends on the  initial condition.
If the latter lies in  the basin of attraction of a certain equilibrium point, the final  configuration will be the one at this specific steady state.
Therefore, it is important to assess, for each possible attractor, the domain of attraction.

In order to have a graphical representation of the separatrix manifolds we study the problem described above in dynamical systems composed by two or three equations. We have analyzed the problem of the reconstruction of the domains of attraction for competition models of dimension two and three which present bistability in \cite{Cavoretto11,C-D-P-V}.

In this paper we instead focus on systems of two or three ordinary differential equations aimed at finding the basins of attraction of three different stable equilibria. Moreover, we propose the problem of the separatrix curves or surfaces that determine the domains of attraction, when bistability occurs, as a particular case of the algorithm described in this article. 
Thus, the  aim of this work is to construct  approximation  surfaces which partition the phase space into three regions.

The approximation of the attraction basins leads to a method consisting of two steps:
\begin{enumerate}
	\item detection of the points lying on the separatrix manifolds,
	\item interpolation of such points in a suitable way, \cite{carr97,Chen,Cuomo,hoppe94,hoppe92,Turk02}.
\end{enumerate}
For this purpose we have implemented several \textsc{Matlab} routines  for the approximation of the points lying on the manifolds determining the basins of attraction, obtained by a bisection algorithm, and for the graphical representation of such manifolds. The \textsc{Matlab} software here discussed is available at:
\begin{center}
	\emph{http://hdl.handle.net/2318/1520518}.
\end{center}

The separatrix manifolds generated by a saddle point are determined locally (by linearization) in well-known examples, see e.g.  \cite{Hale}.
Specifically, even if some techniques to prove the existence of invariant sets have already been developed, none of them, except for particular and well-known cases, allows to have a graphical representation of the separatrix manifolds, \cite{Dellnitz,Johnson}.
Such  techniques are based on results from algebraic topology, and thus such methods are not constructive in the sense that they do not give a precise structure and location of the invariant sets.  Furthermore, numerical tools based on characterizing in (exponentially) asymptotically autonomous systems  a Lyapunov function as a solution of a suitable linear first-order PDE have already been developed. Such  equation is then approximated using meshless collocation methods,  \cite{Giesl11,Giesl12}.

Our aim is instead more ambitious since, on the contrary, the software presented here allows to  reconstruct the basin of attraction of each equilibrium in a three-dimensional dynamical system, providing a graphical representation of the separatrix curves or surfaces. Moreover,  we are not restricted to asymptotically autonomous systems and thus the transformations made in order to use powerful methods, which are well-suited only for autonomous models, are not here necessary.

The paper is organized as follows. In Section \ref{pum} we describe the method used for
approximating the manifolds determining the basins of attraction.
Section \ref{det}  is devoted to the presentation of the designed algorithms for the detection of
points lying on such surfaces.  Section \ref{ne} contains our numerical results.

\section{Approximation of manifolds determining the basins of attraction: interpolation phase}
\label{pum}

In this section we describe the method used for the reconstruction of the basins of attraction. The latter are often described by implicit equations, consequently we use an implicit scheme, specifically the implicit partition of unity method, to reconstruct the domains of attraction. 

We will describe such method for a 3D dataset, but it can easily be adapted to a 2D dataset so as to allow the reconstruction also of implicit curves \cite{Fasshauer}. However, an implicit approach is not always necessary for the approximation of some separatrix manifolds. In fact we have already obtained good results with the explicit partition of unity method, \cite{Cavoretto11,C-D-P-V}. But the nature of a curve or of a surface is known only after detecting the points lying on them. Thus, though such curves or surfaces might usually be expressed by an explicit equation, we use the more general implicit partition of unity technique. The use of such method is the key step which allows to reconstruct the basins of attraction of any stable equilibrium point.

\subsection{Implicit surface reconstruction}\label{PUMi}

Given a point cloud data set, i.e. data in the form $ {\cal X}_N= \{ \boldsymbol{x}_i \in \mathbb{R}^{3},$ $ i=1, \ldots, N \}$, belonging to an unknown two dimensional manifold $\mathscr{M}$, namely a surface in $ \mathbb{R}^{3}$,  we seek another
surface $ \mathscr{M}^{*}$ that is a reasonable approximation to $\mathscr{M}$. For the implicit approach, we think of $\mathscr{M}$
as the surface of all points $\boldsymbol{x}  \in \mathbb{R}^{3}$ satisfying the implicit equation:
\begin{equation}
\label{eqimp}
f(\boldsymbol{x} )=0,
\end{equation}
for some function $f$, which implicitly defines the surface $\mathscr{M}$ \cite{Fasshauer}. This means that the equation \eqref{eqimp} is the zero iso-surface of the trivariate function $f$, and therefore this iso-surface coincides with $\mathscr{M}$.
The surface $\mathscr{M}^{*}$ is constructed via partition of unity interpolation, \cite{Fasshauer},   as shown in Subsection \ref{PUMe}. Unfortunately, the solution of this problem, by imposing the interpolation conditions \eqref{eqimp}, leads to the trivial solution, given by the identically zero function, \cite{Cuomo}.
The key to finding the interpolant of the trivariate function $f$, from the given data points $ \boldsymbol{x}_i,$ $ i=1, \ldots, N,$ is to use additional significant interpolation conditions, i.e. to add an extra set of \emph{off-surface points}. Once we define the augmented data set, we can then compute a three dimensional interpolant $\cal{I}$ to the total set of points, \cite{Fasshauer}.
A common practice is to assume that
in addition to the point cloud data
the set of surface oriented normals $\boldsymbol{n}_i \in \mathbb{R}^{3} $ to the surface
$ \mathscr{M}$  at the points  $ \boldsymbol{x}_i$ is also given.
We construct the extra off-surface points by
taking a small step away along the surface normals, i.e. we obtain for each data point
$ \boldsymbol{x}_i$ two additional off-surface points. One point lies \emph{outside} the manifold $ \mathscr{M}$  and is
given by
$$
\boldsymbol{x}_{N+i}=\boldsymbol{x}_i+ \delta \boldsymbol{n}_i,
$$
whereas the other point lies \emph{inside} $ \mathscr{M}$ and is given by
$$
\boldsymbol{x}_{2N+i}= \boldsymbol{x}_i- \delta \boldsymbol{n}_i,
$$
$\delta$ being the stepsize. The union of the sets $\cal{X}_{ \delta}^{+}=$ $  \{ \boldsymbol{x}_{N+1},$ $\ldots, \boldsymbol{x}_{2N} \}$, $\cal{X}_{ \delta}^{-}=$ $ \{ \boldsymbol{x}_{2N+1},\ldots,$ $\boldsymbol{x}_{3N} \}$ and ${\cal X}_N$  gives the overall set of points on
which the interpolation conditions are assigned. Note that if we have zero normals in the given normal data set, we must exclude such points, \cite{Fasshauer}.

Now, after creating the data set, we compute the interpolant $\cal{I}$
whose zero contour (iso-surface $\cal{I}=$ $0$) interpolates the given point cloud data, and whose iso-surfaces $\cal{I}=$ $1$ and $\cal{I}=$ $-1$ interpolate 
$\cal{X}_{ \delta}^{+}$ and $\cal{X}_{ \delta}^{-}$.
The values $+1$ or $-1$ are arbitrary. Their precise value is not as critical as the choice of $\delta$. In fact the stepsize can be rather critical for a good surface fit, \cite{carr01,Fasshauer}. A suitable value for such parameter will be discussed in Section \ref{ne}. Finally, we just render the resulting approximating surface $ \mathscr{M}^{*}$ as the zero contour of the 3D interpolant, \cite{Fasshauer}. If the normals are not explicitly given, see \cite{hoppe94,hoppe92}. 

\subsection{Partition of unity method and radial basis function interpolation} 
\label{PUMe}
In Subsection \ref{PUMi} we have presented an approach, to obtain a surface that fits the given 3D scattered data set,  based on the use of implicit surfaces defined in terms of some meshfree approximation methods such as the partition of unity interpolation,  \cite{Cavoretto13,Chen,Fasshauer,Melenk96,Wendland02a,Wendland05}. 

Let ${\cal X}_N=\{ \boldsymbol{x}_i \in \mathbb{R}^{3}, i = 1, \ldots , N \}$ be a set of distinct data points or nodes, arbitrarily distributed in a domain $\Omega \subseteq \RR^3$, with an associated set ${\cal F}_N=\{f_i, i = 1,  \ldots , N \}$ of data values or function values, which are obtained by sampling some (unknown) function $f:\Omega \rightarrow \RR$ at the nodes, i.e., $f_i=f(\boldsymbol{x}_i)$, $i=1,\ldots,N$. 

The basic idea of the partition of unity method is to start with a partition of the open and bounded domain $\Omega$ into $d$  subdomains $\Omega_j$ such that $\Omega \subseteq \bigcup_{j=1}^{d} \Omega_j$ with some mild overlap among the subdomains. 

Associated with these subdomains we choose a partition of unity, i.e. a family of compactly supported, non-negative, continuous functions $W_j$ with $\text{supp}(W_j) \subseteq \Omega_j$ such that 
\begin{equation}
\sum_{j=1}^{d} W_j(\boldsymbol{x}) = 1, \hspace{1cm} \boldsymbol{x} \in \Omega.
\end{equation}
The global approximant thus assumes the following form
\begin{eqnarray}
\label{pui}
{\cal I}(\boldsymbol{x})= \sum_{j=1}^{d} R_j(\boldsymbol{x}) W_j(\boldsymbol{x}), \hspace{1cm} \boldsymbol{x} \in \Omega.
\end{eqnarray}
For each subdomain $\Omega_j$ we define a local \textsl{radial basis function} interpolant \cite{Iske11} $R_j:\Omega \rightarrow \RR$ of the form 
\begin{eqnarray}
\label{intfun}
R_j({\boldsymbol{x}})=\sum_{k=1}^{N_j} c_{k} \phi (||\boldsymbol{x}-\boldsymbol{x}_k||_2),
\end{eqnarray}
where $\phi:[0,\infty) \rightarrow \RR$ is called \textsl{radial basis function}, $||\cdot||_2$ denotes the Euclidean norm, and $N_j$ indicates the number of data points in $\Omega_j$. Moreover, $R_j$ satisfies the interpolation conditions
\begin{eqnarray}
\label{condinterp} 
R_j(\boldsymbol{x}_i)=f_i, \hspace{1.cm} i=1,\ldots,N_j.
\end{eqnarray} 
In particular, we observe that if the local approximants satisfy the interpolation conditions \eqref{condinterp}, then the global approximant also interpolates at $\boldsymbol{x}_i$, i.e. ${\cal I}(\boldsymbol{x}_i)=f(\boldsymbol{x}_i)$, for $i=1,\ldots,N_j$.  

Solving the $j$-th interpolation problem \eqref{condinterp} leads to a system of linear equations of the form
\begin{align*}
\left[
\begin{array}{cccc}
\phi (||\boldsymbol{x}_1-\boldsymbol{x}_1||_2)   & \cdots  & \phi (||\boldsymbol{x}_1-\boldsymbol{x}_{N_j}||_2)        \\
\vdots         & \vdots           & \vdots        \\
\phi (||\boldsymbol{x}_{N_j}-\boldsymbol{x}_1||_2)   & \cdots  & \phi (||\boldsymbol{x}_{N_j}-\boldsymbol{x}_{N_j}||_2)        
\end{array}
\right]
\left[
\begin{array}{c}
c_1\\
\vdots\\
c_{N_j}
\end{array}
\right]
=
\left[
\begin{array}{c}
f_1\\
\vdots\\
f_{N_j}
\end{array}
\right],
\end{align*}
or simply
\begin{equation} \label{mat}
\Phi\boldsymbol{c}=\boldsymbol{f}. \nonumber
\end{equation}
Now, we give the following definition (see \cite{Wendland05}).
\begin{definition}
	Let $\Omega \subseteq \RR^3$ be a bounded set. Let $\{\Omega\}_{j=1}^{d}$ be an open and bounded covering of $\Omega$. This means that all $\Omega_j$ are open and bounded and that $\Omega$ is contained in their union. A family of nonnegative functions $\{W_j\}_{j=1}^{d}$ with $W_j \in C^k(\RR^3)$ is called a $k$-stable partition of unity with respect to the covering $\{\Omega_j\}_{j=1}^{d}$ if
	\begin{enumerate}
		\item[1)] ${\rm supp}(W_j) \subseteq \Omega_j$;
		\item[2)] $\sum_{j=1}^{d} W_j(\boldsymbol{x}) \equiv 1$ on $\Omega$;
		\item[3)] for every $\beta \in \NN_0^3$ with $|\beta| \leq k$ there exists a constant $C_{\beta} > 0$ such that
		\begin{equation}
		||D^{\beta}W_j||_{L_{\infty}(\Omega_j)}\leq C_{\beta}/\delta_j^{|\beta|}, \hspace{1.cm} j=1,\ldots,d, \nonumber
		\end{equation}
		where $\delta_j = {\rm diam}(\Omega_j)=\sup_{\boldsymbol{x},\boldsymbol{y} \in \Omega_{j}} ||\boldsymbol{x}-\boldsymbol{y}||_2$.
	\end{enumerate}
\end{definition}

In agreement with the statements in \cite{Wendland02a} we require some additional regularity assumptions on the \textsl{covering} $\{\Omega_j\}_{j=1}^{d}$.

\begin{definition} \label{defpr}
	Suppose that $\Omega \subseteq  \RR^3$ is bounded and ${\cal X}_N=\left\{\boldsymbol{x}_i, i=1,\ldots,N\right\} \subseteq \Omega$ are given. An open and bounded covering $\{\Omega_j\}_{j=1}^{d}$ is called regular for $(\Omega,{\cal X}_N)$ if the following properties are satisfied:
	\begin{itemize}
		\item[(a)] for each $\boldsymbol{x} \in \Omega$, the number of subdomains $\Omega_j$ with $\boldsymbol{x} \in \Omega_j$ is bounded by a global constant $K$;
		\item[(b)] each subdomain $\Omega_j$ satisfies an interior cone condition \cite{Wendland05};
		\item[(c)] the local fill distances $h_{{\cal X}_{N_j}, \Omega_j}$, where ${\cal X}_{N_j}={\cal X}_N \cap \Omega_j$, are uniformly bounded by the global fill distance $h_{{\cal X}_N, \Omega}$, i.e.	
		\begin{eqnarray} 
		h_{{\cal X}_{N}, \Omega} = \sup_{\boldsymbol{x} \in \Omega}\min_{\boldsymbol{x}_k\in {\cal X}_N} ||\boldsymbol{x}-\boldsymbol{x}_k||_2. \nonumber
		\end{eqnarray}	
	\end{itemize}
\end{definition}

Let $C_{\nu}^k(\RR^3)$ be the space of all functions $f \in C^k$ whose derivatives of order $|\beta|=k$ satisfy $D^{\beta}f(\boldsymbol{x})= {\cal O}(||\boldsymbol{x}||_2^{\nu})$ for $||\boldsymbol{x}||_2 \rightarrow 0$. The following convergence result is well known (see, e.g., \cite{Fasshauer,Wendland05}).

\begin{theorem}
	Let $\Omega \subseteq  \RR^3$ be open and bounded and assume that ${\cal X}_N = \{\boldsymbol{x}_i, i=1,$ $\ldots,N \}\subseteq \Omega$. Let $\phi \in C_{\nu}^k(\RR^3)$ be a strictly positive definite function. Let $\{\Omega_j\}_{j=1}^{d}$ be a regular covering for $(\Omega, {\cal X}_N)$ and let $\{W_j\}_{j=1}^{d}$ be $k$-stable for $\{\Omega_j\}_{j=1}^{d}$. Then the error between $f \in {\cal N}_{\phi}(\Omega)$, where ${\cal N}_{\phi}$ is the native space of $\phi$, and its partition of unity interpolant (\ref{pui}) is bounded by
	\begin{equation}
	|D^{\beta}f(\boldsymbol{x}) - D^{\beta}{\cal I}(\boldsymbol{x})| \leq C h_{{\cal X}_N, \Omega}^{(k+\nu)/2 - |\beta|} |f|_{{\cal N}_{\phi}(\Omega)}, \nonumber
	\end{equation}
	for all $\boldsymbol{x} \in \Omega$ and all $|\beta| \leq k/2$. 
\end{theorem}  

If we compare this result with the global error estimates (see e.g. \cite{Wendland05}), we can see that the partition of unity preserves the local approximation order for the global fit. This means that we can efficiently compute large RBF interpolants by solving small RBF interpolation problems and then glue them together with the global partition of unity $\{W_j\}_{j=1}^{d}$. 

\section{Approximation of manifolds determining the basins of attraction}
\label{det}

This section describes the algorithms  implemented for the detection of points lying on the manifolds delimiting the basins of attraction.
In Subsection \ref{det3} we discuss the problem of partitioning the phase space in three subregions, when the system presents three stable equilibria. Moreover, a final remark explains how the algorithm can be easily adapted in case of bistability.

\subsection{Detection of points determining the basins of attraction of three different equilibria}
\label{det3}

In order to approximate  the basins of attraction, when the system admits three stable equilibria, the general  idea is  to find the points lying on the surfaces determining the domains of attraction and finally to interpolate them with a suitable method. 
The steps of the so-called \emph{detection-interpolation} algorithm are summarized in the {\tt{3D-Detec-Interp Algorithm}}.
At first, we need to consider a set of points as initial conditions, then we take points in pairs and we proceed with a bisection routine to determine a point lying on a surface dividing the domains of attraction,  \cite{C-D-P-V}. The simplest idea, which turns out to be also reliable, consists in considering the set of initial conditions in a cube domain $[0,\gamma]^3$, where $\gamma \in \RR^+$. 

Once we take points in pairs, if the two initial conditions evolve towards the first and  second equilibrium point, respectively, a bisection routine is performed to find a point lying on the surface  delimiting the basin of attraction of both the first and the second stable steady state.
In the same way, if the  trajectories of the two initial data converge to the first and  third  stable equilibrium point,   a point  lying on the surface  determining the domain of attraction of  the first and  third attractor is found.
Similarly,  if the  trajectories of the two initial data converge to the second and third   stable equilibrium point,   a point  lying on the surface  determining the domain of attraction of  the second and   third stable steady state is found,  (see {\fontfamily{pcr} \selectfont Step 5} in the {\tt{3D-Detec-Interp Algorithm}}).  It means that the algorithm finds three sets of points, and the latter, taken in pairs, determine the basins of attraction of the three stable equilibria, (see {\fontfamily{pcr} \selectfont Step 6} in  the {\tt{3D-Detec-Interp Algorithm}}).
 Moreover, in the detection-interpolation algorithm, after considering $n$ equispaced initial conditions on each edge of the
cube $[0,\gamma]^3$,  we construct a grid on the faces of the cube, (see {\fontfamily{pcr} \selectfont Step 2} in  the {\tt{3D-Detec-Interp Algorithm}}):
\begin{equation}
\begin{array}{lll}
 P_{i_1,i_2}^{1}=(x_{i_1},y_{i_2},0) & \quad \textrm{and} \quad P_{i_1,i_2}^{2}=(x_{i_1},y_{i_2},\gamma), & \quad i_1, i_2=1, \ldots, n, \\
P_{i_1,i_2}^{3}=(x_{i_1},0,z_{i_2}) & \quad \textrm{and} \quad P_{i_1,i_2}^{4}=(x_{i_1},\gamma,z_{i_2}), & \quad i_1, i_2=1, \ldots, n, \\
P_{i_1,i_2}^{5}=(0,y_{i_1},z_{i_2}) & \quad \textrm{and} \quad P_{i_1,i_2}^{6}=(\gamma,y_{i_1},z_{i_2}), & \quad i_1, i_2=1, \ldots, n,
\end{array}
\label{ic}
\end{equation}
and a bisection routine is applied with initial conditions \eqref{ic}, \cite{C-D-P-V}.

We now analyze the inputs of the detection-interpolation algorithm:
\vskip 0.1 cm
$ n \in \mathbb{N}^{+}$: number of equispaced points on each edge of the cube; it is used to define the set of initial data.
\vskip 0.1 cm
$ \gamma \in \mathbb{R}^{+}$:    edge length of the cube.
\vskip 0.1 cm
 $tol \in \mathbb{R}^{+}$:  tolerance used during the bisection routine. The latter stops when the distance between the last two generated midpoints is less than $tol$. 
\vskip 0.1 cm
 $t \in \mathbb{R}^{+}$:  integration time, used during an integration routine. We cannot provide a suitable choice for this parameter because it depends on the dynamical system, but the algorithm checks if the allowed integration time is sufficient.
\vskip 0.1 cm
 $\boldsymbol{par}\in \mathbb{R}_{+}^l$:  vector of model parameters, where $l$ is the number of model  parameters. The latter must be chosen so that the system presents three stable equilibria.
\vskip 0.1 cm
  $npi \in \mathbb{N}^+$:   number of equilibrium points to be interpolated, typically the origin when it is unstable and the point having all non zero coordinates, such as in population dynamics the coexistence equilibrium, if it is a saddle point.\footnote{In case of two equilibria a saddle point partitions the phase space into two regions, called the basins of attraction of the equilibria. In case of three equilibria instead, several saddles are involved in the dynamics. But the three separating surfaces intersect together at only one saddle which corresponds to a point where all the populations are nonnegative.
}
\vskip 0.1 cm
 $E \in \mathbb{R}^{M \times 3}$:     matrix of the equilibria, where $M$ is the total number of equilibria. 


\vskip 0.1 cm
$\boldsymbol{\varepsilon} \in \mathbb{R}_+^3$: vector of  shape parameters for RBFs used for the reconstruction, via interpolation, of the three surfaces determining the basins of attraction, (see Subsection \ref{PUMe}). Thus it is an input of the interpolation routine, (see {\fontfamily{pcr} \selectfont Step 7} in the {\tt{3D-Detec-Interp Algorithm}}). 
\vskip 0.1 cm
    $\boldsymbol{d}^{PU} \in \mathbb{N}_+^3$:   the number  of subdomain points along one direction of  $\Omega$. 
  It is an input of the interpolation routine,  (see {\fontfamily{pcr} \selectfont Step 7} in the {\tt{3D-Detec-Interp Algorithm}}), which is used to construct  a grid of equally spaced centres of partition of unity subdomains, \cite{Cavoretto14a}. Therefore, $d_i$ can be chosen taking  $ d_i= (d_i^{PU})^3$, $i=1,2,3$, where $d_i$ is the number of subdomains points used to approximate the $i$-th surface.
\vskip 0.1 cm
 $\boldsymbol{K} \in \mathbb{N}_+^3$: vector containing the number of the nearest points used to estimate the normals; input of the interpolation routine,  (see {\fontfamily{pcr} \selectfont Step 7} in  the {\tt{3D-Detec-Interp Algorithm}}). 
\vskip 0.1cm
\begin{table}[!t]
\captionsetup{labelformat=empty}
\begin{center}
\begin{tabular}{p{11.5cm}*{1}{c}}
\hline
\vskip 0.01 cm 
{\fontfamily{pcr} \selectfont Step 1:}
Check if the system presents exactly three stable equilibria.
\vskip 0.12 cm
{\fontfamily{pcr} \selectfont Step 2:}
 Definition of initial conditions on the faces of the cube of edge $\gamma$.
\vskip 0.08 cm
\hskip 2.2 cm Set 
 $s=1$.
\vskip 0.12 cm
{\fontfamily{pcr} \selectfont Step 3:} \textbf{While} $s<=5$
\vskip 0.12 cm
\hskip 0.55 cm {\fontfamily{pcr} \selectfont Step 4:} \textbf{For} $i_1=1, \ldots, n$
\vskip 0.12 cm
\hskip 1.2 cm {\fontfamily{pcr} \selectfont Step 5:} \textbf{For} $i_2=1, \ldots, n$
\vskip 0.08 cm
\hskip 3.3 cm $[\boldsymbol{q}^3,\boldsymbol{q}^2,\boldsymbol{q}^1]=$ BISECTION$(P_{{i_1},{i_2}}^s,P_{{i_1},{i_2}}^{s+1},t,tol,\boldsymbol{par},E )$, 
\vskip 0.08 cm
\hskip 3.3 cm the matrices $Q^3_{j_3,1}$, $Q^2_{j_2,1}$, $Q^1_{j_1,1}$ are uploaded.
\vskip 0.12 cm
\hskip 2.2 cm 
$s=s+2$.
\vskip 0.12 cm
{\fontfamily{pcr} \selectfont Step 6:} Define $Q^1_{J_1,k}=[Q^2_{j_2,k};Q^3_{j_3,k}]$, $Q^2_{J_2,k}=[Q^1_{j_1,k};Q^3_{j_3,k}]$, 
\vskip 0.08 cm
 \hskip 2.2 cm $Q^3_{J_3,k}=[Q^1_{j_1,k};Q^2_{j_2,k}]$,  $k=1,2,3$.
\vskip 0.12 cm
{\fontfamily{pcr} \selectfont Step 7:} ${\cal I}_1=$ INTERPOLATION$({\cal Q }^1_{J_1,k},\varepsilon_1,d_1^{PU},K_1)$,
\vskip 0.08 cm
\hskip 1.65 cm  ${\cal I}_2=$ INTERPOLATION$({\cal Q }^2_{J_2,k},\varepsilon_2,d_2^{PU},K_2)$, 
\vskip 0.08 cm
 \hskip 1.65 cm ${\cal I}_3=$ INTERPOLATION$({\cal Q }^3_{J_3,k},\varepsilon_3,d^{PU}_3,K_3)$,  where ${\cal Q }^1_{J_1,k},$ 
\vskip 0.08 cm
\hskip 2.2 cm ${\cal Q }^2_{J_2,k},$   ${\cal Q }^3_{J_3,k}$ are the points found by the detection 
\vskip 0.08 cm
\hskip 2.2 cm  algorithm and  the $npi$  equilibrium points to be interpolated. 
\\[\smallskipamount]
\hline
\end{tabular}
\end{center}
\caption{The {\tt{3D-Detec-Interp Algorithm.}} The detection-interpolation pseudo-code. It summarizes the steps needed to determine the points lying on the surfaces delimiting the  basins of attraction and to reconstruct these
surfaces.}
\label{algoritmo1}
\end{table}

More in detail, after defining the initial data \eqref{ic}, a bisection-like  routine is performed with the latter,  (see {\fontfamily{pcr} \selectfont Step 2-5}  in  the {\tt{3D-Detec-Interp Algorithm}}). Such routine integrates the system with a pair of initial conditions for the parameter set $\boldsymbol{par}$ in an time interval $t$. Then it checks, among the equilibria stored in the matrix $E$, where the trajectories originating in such initial conditions ultimately stabilize. 
Then if the two initial conditions evolve toward two different stable equilibria it  provides a point named:
\begin{enumerate}
\item
$\boldsymbol{q}^3$, if the point lies   on the surface delimiting the domain of attraction of both the first and  the second equilibrium point, or
\item
 $\boldsymbol{q}^2$, if lies   on the surface determining the basin of attraction of both the first and  the third equilibrium point, or
\item $\boldsymbol{q}^1$, if the point lies   on the surface delimiting the domain of attraction of both the second and  the third equilibrium point.
\end{enumerate}

Summarizing, once we apply the bisection algorithm with initial conditions \eqref{ic}, three different sets of points are detected.
These sets,  in pairs, identify the three basins of attraction. Considering then the method described in Section \ref{pum}, the associated algorithm interpolates such points and returns values of the interpolants ${\cal I}_1$, ${\cal I}_2$ and  ${\cal I}_3$. They approximate the basins of attraction of the first, second and third equilibrium point, respectively, (see {\fontfamily{pcr} \selectfont Step 7} in the  {\tt{3D-Detec-Interp Algorithm}}).

\begin{remark}
The separatrix surfaces in case of bistability, investigated in \cite{C-D-P-V}, can be obtained as a particular  case of the detection-interpolation algorithm analyzed in this section. Since in this case we obtain only one manifold, the input parameters $\varepsilon$, $d^{PU}$, $K$ in the detection-interpolation algorithm, are scalar values.
\end{remark}

\begin{remark}
	The approximation of the basins of attraction  for two dimensional dynamical systems easily follows from the detection-interpolation algorithm.
	In the 2D case, we start considering $n$ equispaced initial conditions on each edge of the
	square $[0,\gamma]^2$; thus the bisection routine is applied with the following initial conditions, \cite{C-D-P-V}:
	\begin{equation}
	\begin{array}{lll}
	P_{i}^{1}=(x_i,0)  & \quad \textrm{and} \quad P_{i}^{2}=(x_i,\gamma), & \quad i=1, \ldots, n, \\
	P_{i}^{3}=(0,y_i)  & \quad \textrm{and} \quad P_{i}^{4}=(\gamma,y_i), & \quad i=1, \ldots, n. 
	\end{array}
	\label{2ic}
	\end{equation}
\end{remark}

\section{Numerical experiments}
\label{ne}


In this section we summarize the extensive experiments performed to test our detection and approximation techniques. Specifically, in Subsection \ref{ne3} and \ref{ne2} respectively, we test the routines for  3D and 2D dynamical systems, considering the cases in which such models admit both two and three stable equilibria.

For the dynamical systems in consideration we establish conditions to be imposed on the parameters so that  the separatrix manifolds exist. Here, after detecting the points lying on the latter with the algorithm described in Section \ref{det}, at first we  compute the normal vectors and consistently orient them to the surfaces by choosing, for the three different manifolds, the  nearest neighbours $K_i$, $i=1,2,3$. Typically we set   $K_i$, $i=1,2,3$ between $5$ and $10$.
Then we  build  the extra off-surface points by marching a small distance $ \delta$ along the surface normals, as shown in Subection \ref{PUMi}; following \cite{Fasshauer}, we take $\delta $ to be $1 \%$  of the maximum dimension of the data. Finally we interpolate the points lying on the separatrix surfaces with the implicit partition of unity  method,   described in Subsection \ref{PUMe}, using in \eqref{intfun} the compactly supported Wendland's $C^{2}$ function, \cite{Wendland05}, as local approximants 
\begin{displaymath}
\phi(r)=(1-\varepsilon r)^{4}_{+} (4 \varepsilon r+1).
\end{displaymath}
Here $r = \left\|\cdot\right\|_2$ is the Euclidean norm, $(\cdot)_+$ denotes the truncated power function and $\epsilon > 0$ is the shape parameter. Such parameter determines the size of the support of the basis function. Its choice can significantly affect the final result. Specifically, for the three manifolds, we choose the shape parameters so that $0.01 \leq \varepsilon_i \leq 0.1$, $i=1, 2, 3$.
Assuming to have a nearly uniform node distribution such as the Halton points, according to \cite{Cavoretto14a},
a possible choice for the number of subdomains centers consists in constructing a uniform grid of $d_i=(d^{PU}_i)^s$ centers, where $s$ is the dimension of the dynamical system and  $d^{PU}_i= \lceil{ 1/2 (N_i/2)^{ \frac{1}{s}}} \rceil$, $i=1, 2,3$. However, we point out that in our tests we find good results even with different $d^{PU}_i$, $i=1,2,3$. This is due to the fact that we deal with concrete and unstructured data. 

Such choices,  described above, are suitable assuming to start  with $ 8 \leq n \leq 15$ equispaced initial conditions on each edge of the cube $[0,\gamma]^3$.
For the tolerance used in the bisection routine, a recommended value is $10^{-3} \leq tol \leq 10^{-5}$, since it allows to achieve a good trade-off between accuracy and computational cost.

\subsection{3D detection-interpolation tests}
\label{ne3}
A model chosen to test the detection-interpolation algorithm is the classical three-populations competition model.
Letting $x$, $y$ and $z$ denote the populations, we consider the following system 
\begin{equation} \label{model3d}
\begin{array}{ll}
\frac{ \displaystyle  dx}{ \displaystyle  dt}=p \big(1- \frac{ \displaystyle  x}{ \displaystyle  u} \big)x-axy-bxz,  & \textrm{} \\
\vspace{.01cm}\\
\frac{ \displaystyle  dy}{ \displaystyle  dt}=q \big(1- \frac{ \displaystyle  y}{ \displaystyle  v} \big)y-cxy-eyz, & \textrm{} \\
\vspace{.01cm}\\
\frac{ \displaystyle  dz}{ \displaystyle  dt}=r \big(1- \frac{ \displaystyle  z}{ \displaystyle  w} \big)z-fxz-gyz,  & \textrm{} 
\end{array}
\end{equation}
where $p$, $q$ and $r$ are the growth rates of $x$, $y$ and $z$, respectively, $a$, $b$, $c$, $e$, $f$ and $g$ are the competition rates,
$u$, $v$ and $w$ are the carrying capacities of the three populations.
The model (\ref{model3d}) describes the interaction of three competing populations within the
same environment (see e.g. \cite{GLMMTV}). 

\begin{figure}[ht!]
	\begin{center}
		\includegraphics[height=.22\textheight]{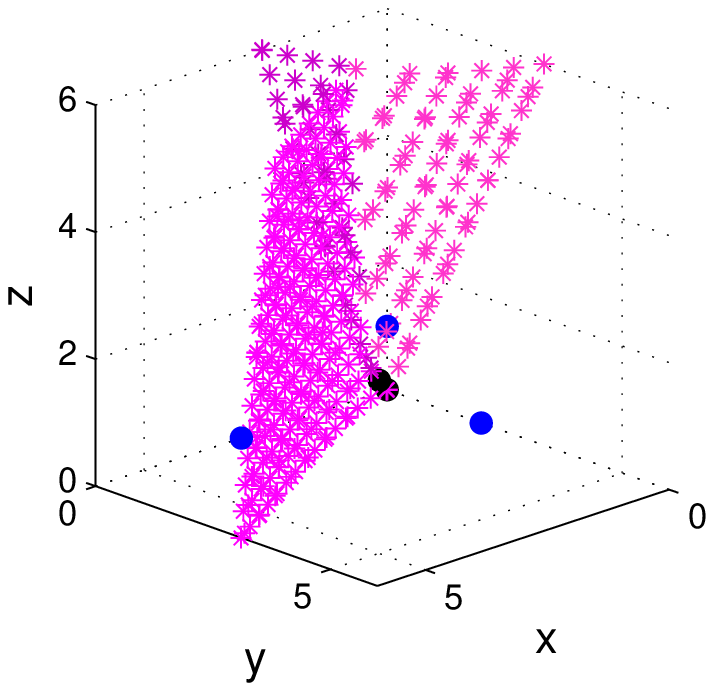} 
		\includegraphics[height=.22\textheight]{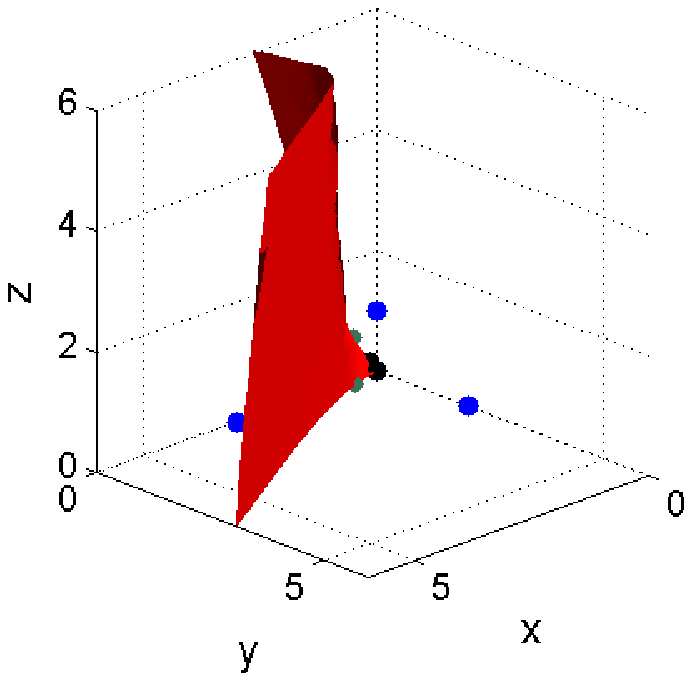} 
		\includegraphics[height=.22\textheight]{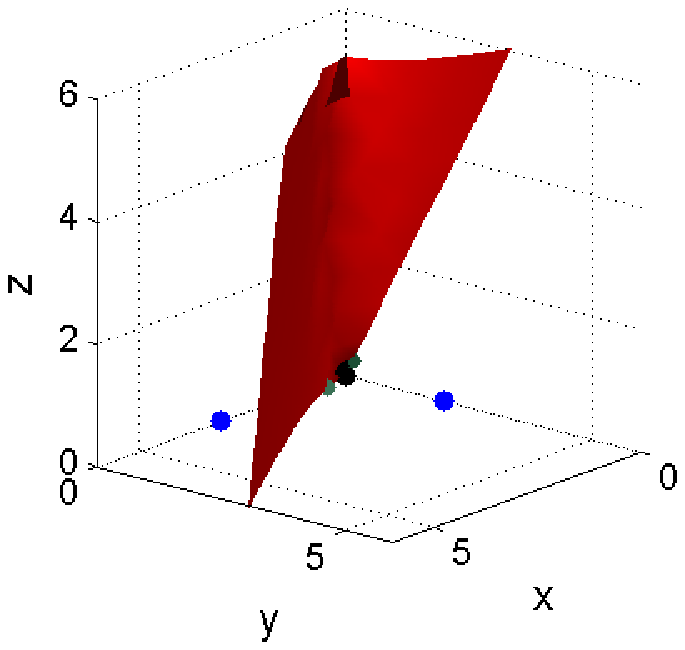} 
		\includegraphics[height=.22\textheight]{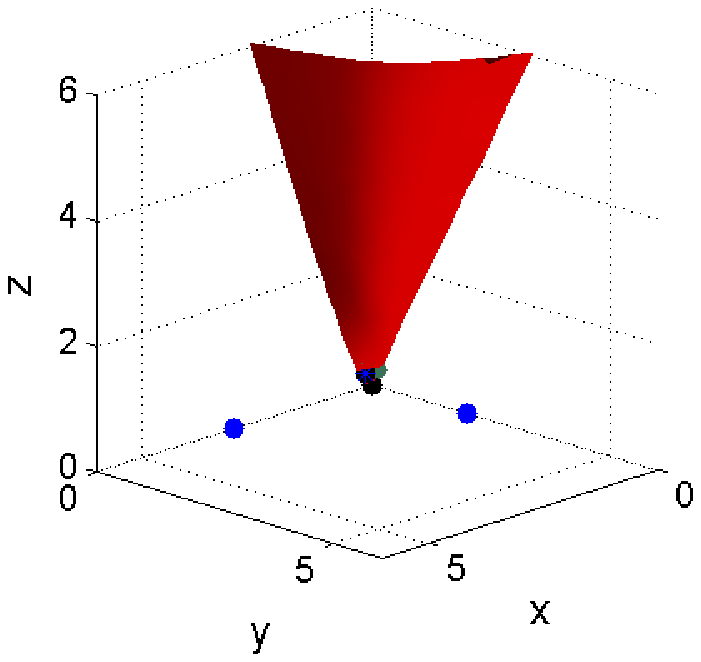} 
		\caption{Set of points lying on the surfaces determining the domains of attraction (top left) and the reconstruction of the basin of attraction of $E_1$ (top right), $E_2$ and $E_3$ (bottom, left to right). The four figures  (left to right, top to bottom)  show the progress of the algorithm: first it generates the points on the separatrices, then in turn each individual basin of attraction. The black and blue circles represent the unstable origin, the coexistence saddle point and the stable equilibria, respectively. Moreover the other saddles ($E_4$, $E_5$ and $E_6$) which, in pairs, lie on the separatrix manifolds of the attraction basins are identified by green circles.
		}
		\label{puntiE12322}
	\end{center}
\end{figure}

\begin{figure}[ht!]
	\begin{center}
		\includegraphics[height=.22\textheight]{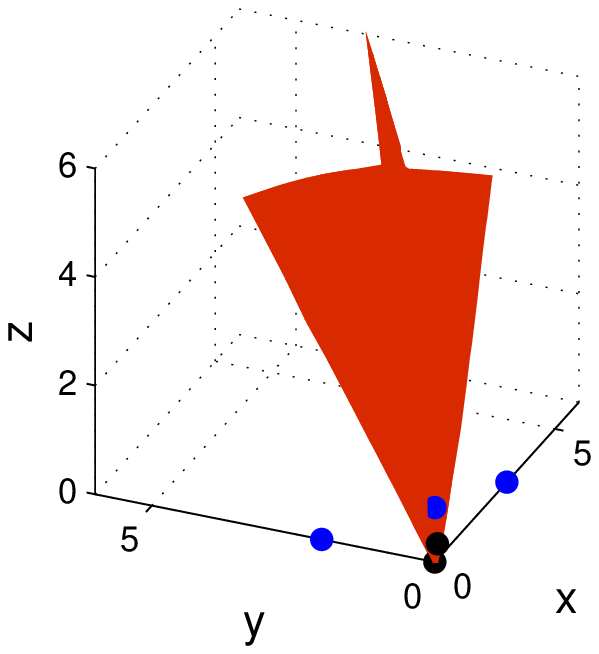} 
		\caption{Reconstruction of the basins of attraction with parameters $p = 1$, $q = 2$, $r =2$, $a = 5$, $b=4$, $c =3$, $e=7$, $f=7$, $g=10$, $u=3$, $v=2$, $w=1$.}
		\label{bacini11}
	\end{center}
\end{figure}

There are eight equilibrium points. The origin $E_0 = (0, 0, 0)$ and the points associated with the survival of only one population $E_1 = (u, 0, 0)$, $E_2 = (0, v, 0)$ and $E_3 = (0, 0, w)$. Then we have the equilibria with two coexisting
populations: 
\begin{displaymath}
\begin{array}{l}
\vspace{0.3em}
E_4 =  \bigg( \frac{\displaystyle uq(av-p)}{\displaystyle cuva-pq},\frac{\displaystyle pv(cu-q)}{\displaystyle cuva-pq},0\bigg), \quad E_5 = \bigg( \frac{\displaystyle ur(bw-p)}{\displaystyle fuwb-rp},0, \frac{\displaystyle wp(fu-r)}{\displaystyle fuwb-rp}\bigg),\\
E_6 = \bigg( 0, \frac{\displaystyle vr(we-q)}{\displaystyle gvwe-qr}, \frac{\displaystyle wq(vg-r)}{\displaystyle gvwe-qr}\bigg).
\end{array}
\end{displaymath}
Finally we have the coexistence equilibrium,
\begin{displaymath}
\left.
\begin{array}{ll}
\vspace{0.3em}
E_7 = & \bigg(\frac{\displaystyle u[p(gvwe-qr)-avr(we-q)-bwq(vg-r)]}{\displaystyle p(gvwe-qr)+uva(rc-fwe)+uwb(fq-gcv)},\\
\vspace{0.3em}
& \frac{ \displaystyle v[q(fuwb-pr)-rcu(wb-p)-pew(fu-r)]}{ \displaystyle q(fuwb-pr)+cuv(ra-gwb)+evw(gp-afu)},\\
&\frac{ \displaystyle r[(cuva-pq)-gpv(cu-q)-ufq(va-p)]}{\displaystyle r(cuva-pq)+bwu(fq-vcg)+evw(gp-fua)} \bigg).
\end{array}
\right.
\end{displaymath}
Letting $p = 1$, $q = 2$, $r =2$, $a = 5$, $b=4$, $c =3$, $e=7$, $f=7$, $g=10$, $u=3$, $v=2$, $w=1$, the  points associated with the survival of only one population, i.e. $E_1 = (3, 0, 0)$, $E_2 = (0, 2, 0)$ and $E_3 = (0, 0, 1)$, are stable, the origin $E_0 = (0, 0, 0)$ is an unstable equilibrium and the coexistence equilibrium $E_7 \approx (0.1899,    0.0270,    0.2005)$ is a saddle point. The remaining equilibria $E_4 \approx (0.6163,0.1591,0)$, $E_5 \approx (0.2195,0,0.5317)$ and $E_6 \approx (0,0.1714,0.2647)$ are other saddle points. The manifolds joining these saddles partition the phase space into the different basins of attraction, but intersect only at the coexistence saddle point, labeled $E_7$.\footnote{In case of bistability the manifold through the origin and a saddle point partitions the phase space into two regions. In case of a  system with three equilibria instead, more saddles are involved in the dynamics. But the three separating manifolds all intersect only at one saddle with all nonnegative populations.}
In this situation we can use the detection-interpolation routine to approximate the basins of attraction.
More precisely, we choose $n=15$, $\gamma=6$, $tol=10^{-3}$, $t=90$,  $\boldsymbol{\varepsilon}=(0.1,0.09,0.08)$, $ \boldsymbol{d}^{PU}=(3,4,4)$,  $\boldsymbol{K}=(7,8,6)$.
Figure \ref{puntiE12322} shows the separatrix points and   the basins of attraction of $E_1$, $E_2$ and $E_3$, (left to right, top to bottom).
Finally, in Figure \ref{bacini11}   we plot together the three basins of attraction.\\

To test our detection-interpolation routine when bistability occurs, we consider the following model, describing a three level food web, with a top predator indicated by $ W$,
the intermediate population $V$ and the bottom prey $N$ that is affected by
an epidemic. It is subdivided into the two subpopulations of susceptibles
$S$ and infected $I$, \cite{derossi},
\begin{equation} \label{model4d}
\begin{array}{ll}
\frac{ \displaystyle  dW}{ \displaystyle  dt}=-mW + pVW,  & \textrm{} \\
\vspace{.01cm}\\
\frac{ \displaystyle  dV}{ \displaystyle  dt}=-lV + eSV - hVW  + qIV, & \textrm{} \\
\vspace{.01cm}\\
\frac{ \displaystyle  dI}{ \displaystyle  dt}= \beta IS - nIV - \gamma I - \nu I,  & \textrm{} \\
\vspace{.01cm}\\
\frac{ \displaystyle  dS}{ \displaystyle  dt}= aS \big(1- \frac{ \displaystyle  S+I}{ \displaystyle  K} \big) - c V S - \beta SI +  \gamma I,  & \textrm{} 
\end{array}
\end{equation}
where $m$ and $l$ are the mortality rates of $W$ and $V$ respectively, $\nu$ is the natural plus disease-related mortality for the bottom
prey, $p$ and $h$ are the predation rates. The disease, spreading by contact at rate $ \beta$, can be overcome, so that infected return to class $S$ at rate $\gamma$. Then the gain obtained by the intermediate population
from hunting of susceptibles is denoted by $e$, which must clearly be smaller
than the damage inflicted to the susceptibles $c$, i.e. $e < c$, the corresponding
loss rate of infected individuals in the lowest trophic level due to capture by
the intermediate population is $n$, while $q < n$ denotes the return obtained by
$V$ from capturing infected prey. In this lowest trophic level, only the healthy prey reproduce at
net rate $a$, while the prey environment carrying capacity is $K$.
The equilibria are the origin $E_0=(0,0,0,0)$, $E_1=(0,0,0,K)$, the disease-free equilibrium with all the trophic levels $E_2$ and the steady state in which only the intermediate population and the bottom healthy prey thrive $E_3$:
\begin{displaymath}
E_2=\bigg( \frac{ \displaystyle  apKe - mecK - apl}{ \displaystyle  ahp}, \frac{m}{p},0, K \frac{ap - cm}{ap} \bigg),  \quad E_3=\bigg(0, \frac{a(Ke- l)}{ecK} ,0,\frac{l}{e}\bigg).
\end{displaymath}
Then we have the point at which just the bottom prey thrives, with endemic disease, $E_4$ and  two equilibria in which the top predators disappear, $E_5$ and $E_6$:
\begin{displaymath}
\begin{array}{l}
\vspace{0.3em}
E_4= \bigg (0,0,  \frac{ \displaystyle a(K \beta \gamma + k \beta \nu-\gamma^2-2 \gamma \nu - \nu^2)}{  \displaystyle  \beta( a \gamma a \nu + K \beta \nu)},\frac{ \displaystyle  \gamma+ \nu}{  \displaystyle \beta} \bigg), \\
E_{5,6}=\bigg(0,\frac{  \displaystyle \beta \hat{S}-\gamma -\nu}{ \displaystyle n},\frac{  \displaystyle l-e \hat{S}}{  \displaystyle q}, \hat{S} \bigg),
\end{array}
\end{displaymath}
where $\hat{S} $ are the roots of $\tilde{A} S^2 + \tilde{B} S +C=0$.

\begin{figure}[ht!]
	\begin{center}
		\includegraphics[height=.22\textheight]{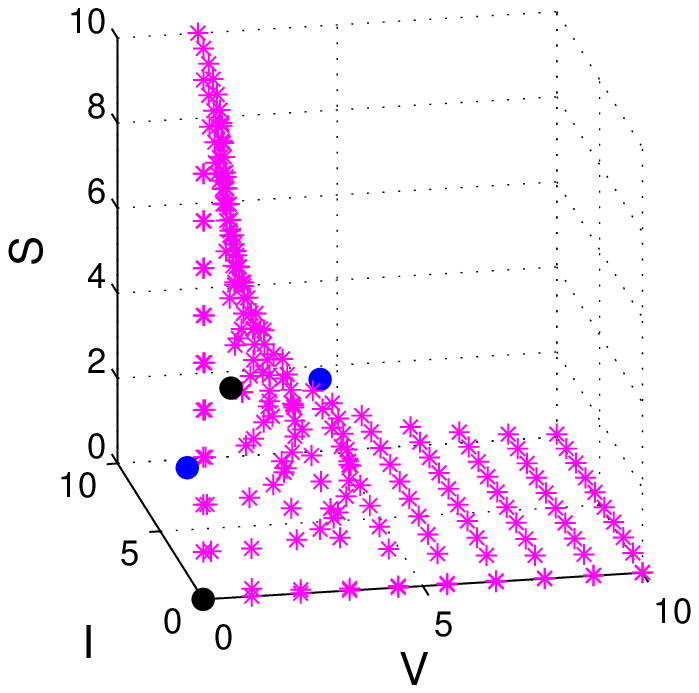} 
		\includegraphics[height=.22\textheight]{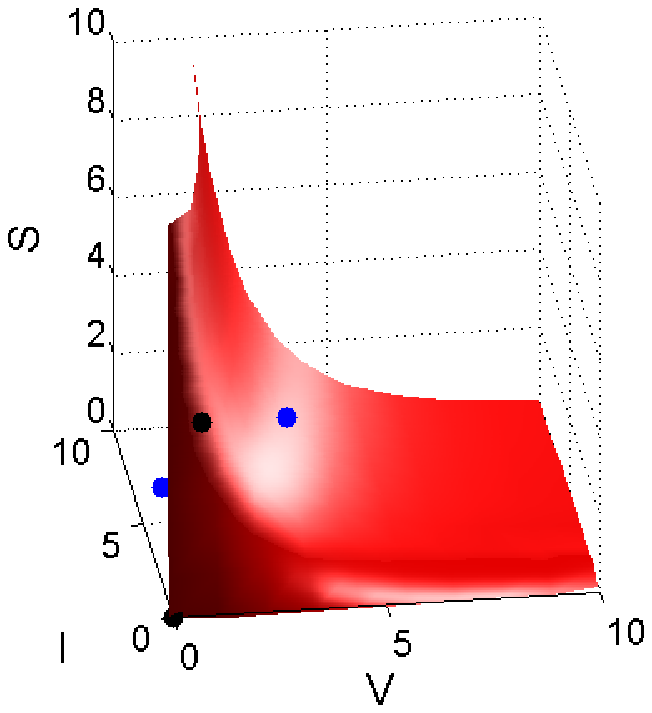} 
		\caption{Set of points lying on  the surface separating the domains of attraction of $E_3$ and $E_4$ (left) and the  reconstruction of the separatrix surface (right). The black and blue circles represent the unstable origin, the saddle point and the stable equilibria,  respectively.
		}
		\label{sep}
	\end{center}
\end{figure}

With the parameters values $l = 10$, $e = 2,$ $q = 1,$ $ \beta = 1.6$, $n = 5$, $\gamma = 1,$ $\nu = 3$, $a = 8$, $K = 6$, $c = 0.5$, the equilibria $E_3 \approx (0,2.6666, 0, 5)$ and $E_4 \approx (0, 0, 1.8421, 2.5)$ are both stable and $E_5 \approx (0, 0.7244, 0.4721, 4.7639)$ is the saddle point that partitions the domain  in the  $W = 0$ three-dimensional phase subspace. Thus system \eqref{model4d} is reduced to a system of three equations and therefore  we can reconstruct the separatrix surface in such subspace with the routine described in Subsection \ref{det3}. The separatrix points and the separatrix surface, shown in Figure  \ref{sep} (left) and  (right) respectively, are the result of the  detection-interpolation  algorithm with $n=11,$ $\gamma=10,$ $tol=10^{-4},$  $t=30$, $\varepsilon=0.6,$ $d^{PU}=4,$ $K=7$.

\subsection{2D detection-interpolation tests}
\label{ne2}
To give an example for a dynamical system of dimension two, we can consider the competition model analyzed in \cite{pastacaldi}. Letting $P$ and $Q$ denote two populations gathering in herds,
we consider the following system describing the competition of  two different populations  within the
same environment:
\begin{equation} \label{model2d}
\begin{array}{ll}
\frac{ \displaystyle  dQ}{ \displaystyle  d \tau}=r \bigg(1- \frac{ \displaystyle  Q}{ \displaystyle  K_Q} \bigg)Q- q \sqrt{  \displaystyle  Q} \sqrt{  \displaystyle P},  & \textrm{} \\
\vspace{.01cm}\\
\frac{ \displaystyle  dP}{ \displaystyle  d \tau}=m \bigg(1- \frac{ \displaystyle  P}{ \displaystyle  K_P} \bigg)P- p \sqrt{  \displaystyle Q} \sqrt{  \displaystyle P}, & \textrm{} 
\end{array}
\end{equation}
where $r$ and $m$  are the growth rates of $Q$ and $P$, respectively, $q$ and $p$ are the competition rates,
$K_Q$,  and $K_P$ are the carrying capacities of the two populations.

Since singularities could arise in the Jacobian when one or both populations
vanish, we define the following new variables,  as suggested in \cite{pastacaldi}:
\begin{equation}\label{tras}
\begin{array}{l}
X(t)= \sqrt{  \displaystyle  \frac{ \displaystyle Q( \tau)}{ \displaystyle K_Q}}, \quad  
Y(t)= \sqrt{  \displaystyle  \frac{ \displaystyle P( \tau)}{ \displaystyle K_P}}, \quad
t= \tau \frac{  \displaystyle  q \sqrt{  \displaystyle K_P}}{  \displaystyle  2 \sqrt{  \displaystyle K_Q}},  \\
\vspace{.01cm}\\
a= \frac{  \displaystyle p  K_Q}{ \displaystyle  q K_P}, \quad
b= \frac{  \displaystyle  r \sqrt{  \displaystyle K_Q}}{  \displaystyle q  \sqrt{  \displaystyle K_P}}, \quad 
c=  \frac{  \displaystyle m  \sqrt{  \displaystyle K_Q}}{  \displaystyle q  \sqrt{  \displaystyle K_P}}.
\end{array}
\end{equation}
Thus the adimensionalized, singularity-free system for  \eqref{model2d} is
\begin{equation} \label{model2dris}
\begin{array}{ll}
\frac{ \displaystyle  dX}{ \displaystyle  d t}=b (1- X^2)X- Y,  & \textrm{} \\
\vspace{.01cm}\\
\frac{ \displaystyle  dY}{ \displaystyle  d t}=c(1-Y^2)Y-aX. & \textrm{} 
\end{array}
\end{equation}
We can  easily verify that the origin $E_0=(0,0)$, and the points associated with the survival of only one population $E_1=(K_Q,0)$ $E_2=(0,K_P)$ are equilibria of \eqref{model2d}. To study the remaining  equilibria we consider  the adimensionalized system, in fact the coexistence equilibria are the roots of the eighth degree equation
$$
cb^3X^8 - 3cb^3X^6 + 3cb^3X^4 - cb(b^2 + 1)X^2 − a + cb = 0.
$$
Observe that in our test we have to take into account that $E_1^{'}=(1,0)$ and $E_2^{'}=(0,1)$, corresponding to $E_1=(K_Q,0)$ and  $E_2=(0,K_P)$ of system \eqref{model2d},  are not critical points of the system \eqref{model2dris}.

\begin{figure}[ht!]
	\begin{center}
		\includegraphics[height=.22\textheight]{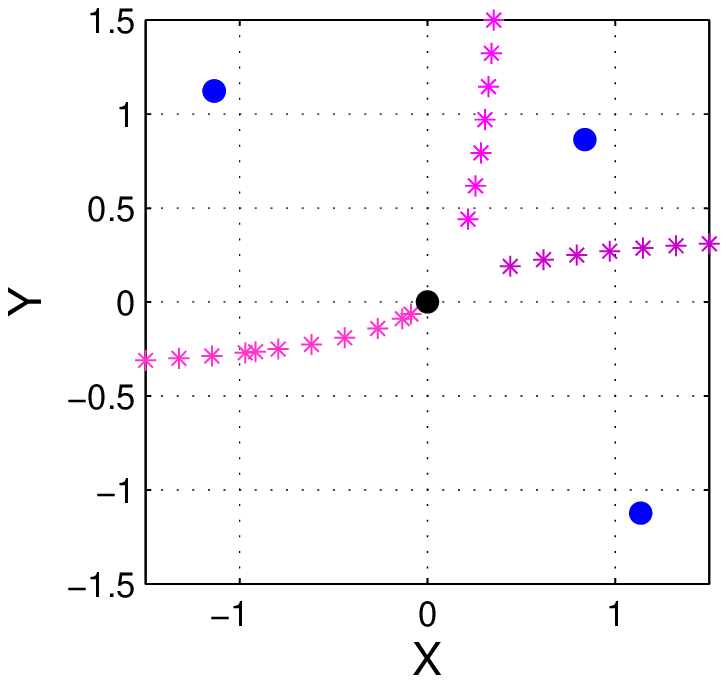} 
		\includegraphics[height=.22\textheight]{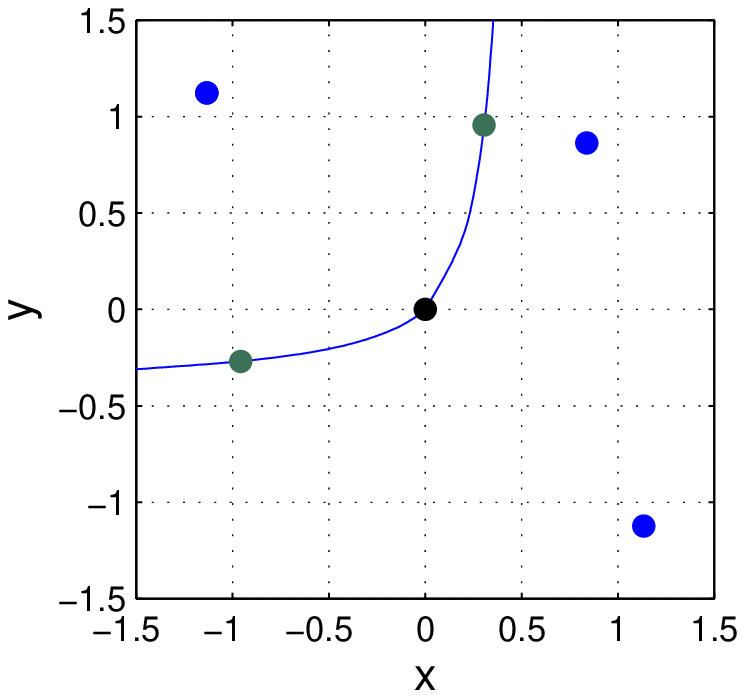} 
		\includegraphics[height=.22\textheight]{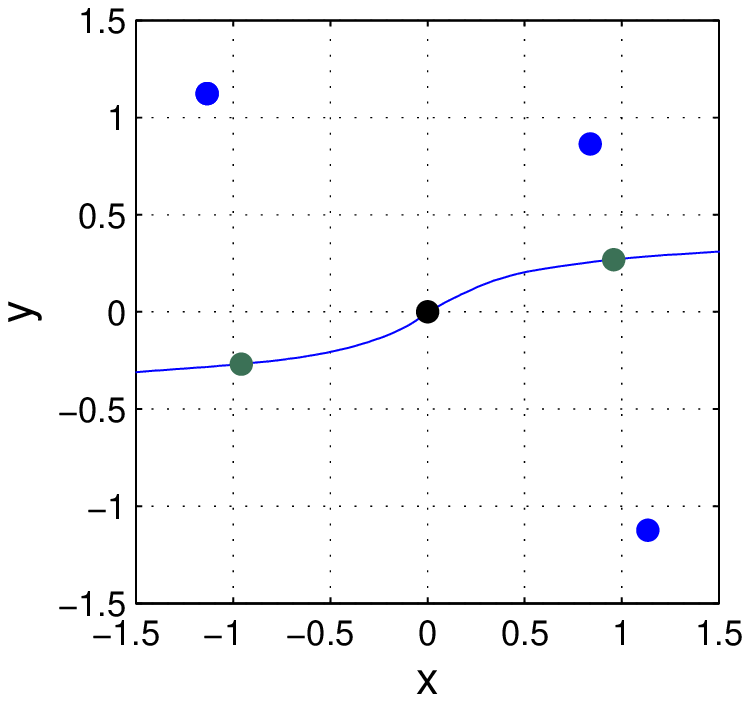} 
		\includegraphics[height=.22\textheight]{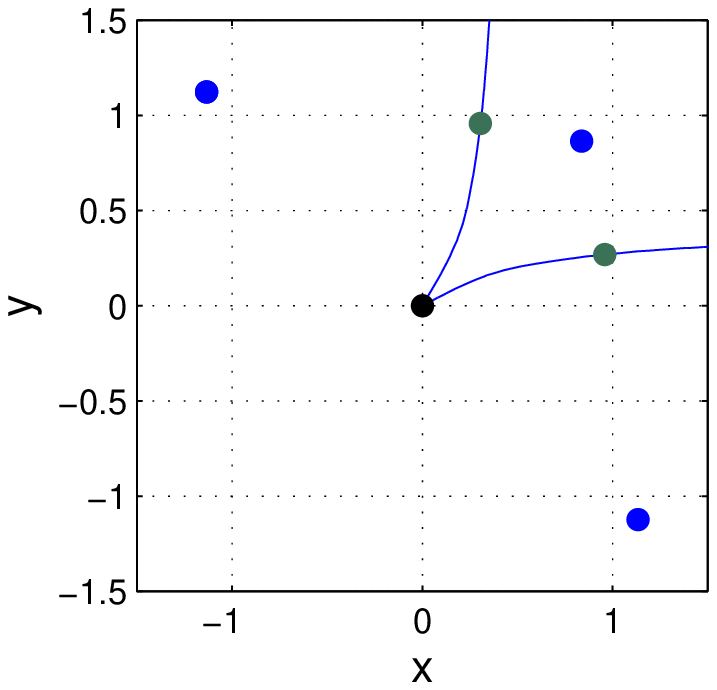} 
		\caption{Set of points lying on the curves determining the domains of attraction (top left) and  the reconstruction of the basin of attraction of $E_1$ (top right), $E_2$ and $E_3$ (bottom, left to right). 
			The four figures (left to right, top to bottom) show the progress of the algorithm: first it
			generates the points on the separatrices, then in turn each individual basin of attraction.
			The black and blue circles represent the  origin and the stable equilibria,  respectively.
			Moreover the other saddles ($E_4,$ $E_5$ and $E_6$) that lie on the separatrix manifolds of the attraction basins are identified by green circles.
		}
		\label{curve_ris_bacini}
	\end{center}
\end{figure}

\begin{figure}[ht!]
	\begin{center}
		\includegraphics[height=.22\textheight]{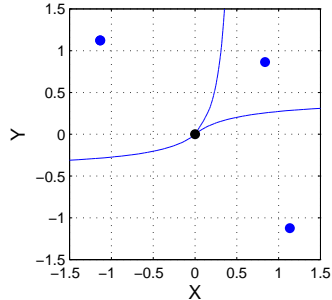} 
		\caption{Reconstruction of the basins of attraction with parameters $r = 0.7895$, $m = 0.7885$, $p = 0.225$, $q = 0.2085$, $Kp = 12$ and $Kq = 10$.
		}
		\label{curve_ris}
	\end{center}
\end{figure}

\begin{figure}[ht!]
	\begin{center}
		\includegraphics[height=.22\textheight]{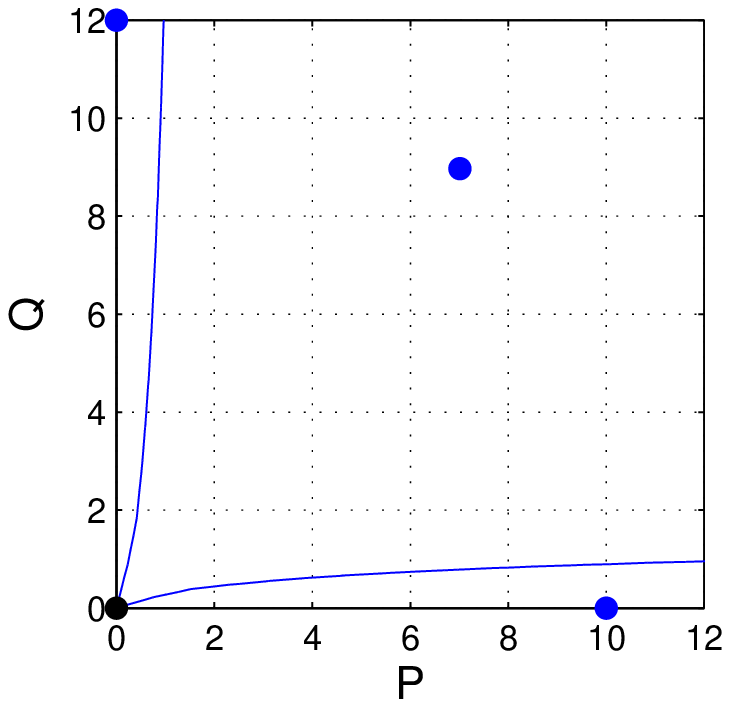} 
		\includegraphics[height=.22\textheight]{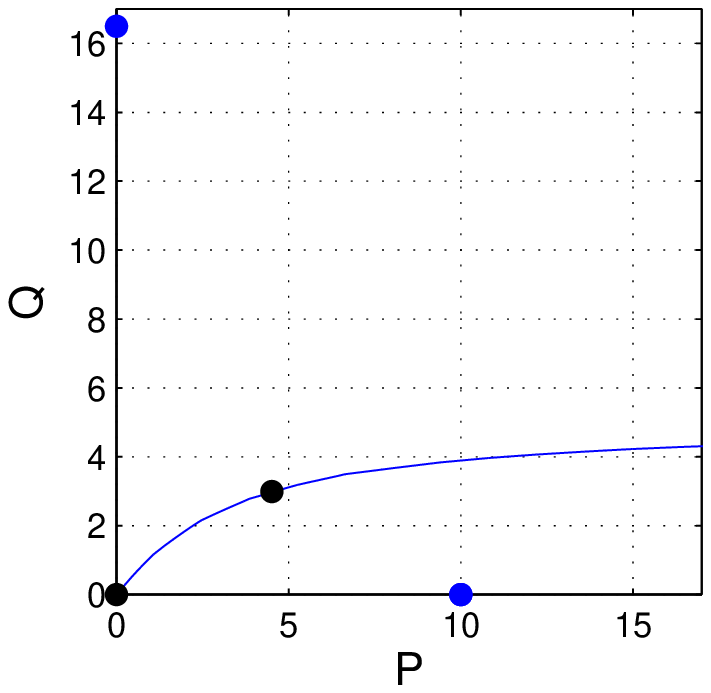} 
		\caption{The basin of attraction of $E_1$, $E_2$ and $E_3$ with parameters $r = 0.7895$, $m = 0.7885$, $p = 0.225$, $q = 0.2085$, $Kp = 12$ and $Kq = 10$ (left), and  the curve separating the basin of attraction of $E_1$, $E_2$ with parameters $r = 0.8888$, $m = 0.602$, $p = 0.401$, $q = 0.5998$, $Kp = 16.5$, $Kq = 10$.
		}
		\label{curve_no_ris}
	\end{center}
\end{figure}

With the parameters $r = 0.7895$, $m = 0.7885$, $p = 0.225$, $q = 0.2085$, $Kp = 12$ and $Kq = 10$,
the points $E_1=(10,0)$, $E_2=(0,12)$
and $E_3 \approx (7.0127, 8.9727)$ are stable equilibria of the system \eqref{model2d}. Instead of integrating  the latter  we consider the  model  \eqref{model2dris}, whose  three stable equilibria  are $E_1^{*} \approx (-1.1342,1.1237)$, $E_2^{*}  \approx (1.1342,-1.1237)$, $E_3^{'} \approx (0.8374,0.8647)$, whereas the origin is the saddle point through which all the three curves go.
Note that when three stable attractors are present there are also other saddles involved in the dynamics, namely $E_4 \approx(-0.9585,-0.2692)$, $E_5 \approx(0.9585,0.2692)$ and $E_6 \approx(0.3055,0.9575)$. Observe that, applying the transformations \eqref{tras}, obviously $E_3^{'}$ corresponds to $E_3$, while $E_1^{*}$ and $E_2^{*}$ are not feasible, but roughly speaking, they represent $E_1$ and $E_2$. In fact the trajectories converging to  $E_1^{*}$ and $E_2^{*}$, under the biological constraint $X \geq 0$, $Y \geq 0$, stop on the axes evolving toward the biological equilibria $E_1^{'}$ and $E_2^{'}$. Therefore we consider $E_1^{*}$, $E_2^{*}$ and $E_3^{'}$.
To apply the algorithm with initial conditions \eqref{2ic} we need a further consideration. Specifically, we have to  translate the problem in the positive plane with the substitutions 
\begin{equation}\label{traslaz}
X^{'}=X+ \frac{ \gamma}{2} \hskip 0.5cm \textrm{and} \hskip 0.5cm  Y^{'}=Y + \frac{ \gamma}{2}, 
\end{equation}
where $ \gamma$  is the length  of the square. At this point we can apply the detection-interpolation algorithm. 
More precisely, we choose:  $n=13$, $\gamma=3$, $tol=10^{-4}$, $t=40$,  $\boldsymbol{\varepsilon}=(0.1,0.06,0.08)$, $ \boldsymbol{d}^{PU}=(4,3,3)$,  $\boldsymbol{K}=(4,6,6)$. Figure \ref{curve_ris_bacini} shows how the algorithm works. It generates first the points lying on the
curves determining the domains of attraction (top left), then subsequently the basins of
attraction of $E_1^{*}$  (top right), $E_2^{*}$ (bottom left) and $E_3^{'}$  (bottom right), in the original  system $X$ and $Y$.
Finally, in Figure  \ref{curve_ris} we plot together the three basins of attraction, always in the original system. Using again the transformation \eqref{tras} we obtain the curves separating the basins of attraction of $E_1$, $E_2$ and $E_3$, shown in Figure \ref{curve_no_ris} (left).

\begin{figure}[h!]
	\begin{center}
		\includegraphics[height=.22\textheight]{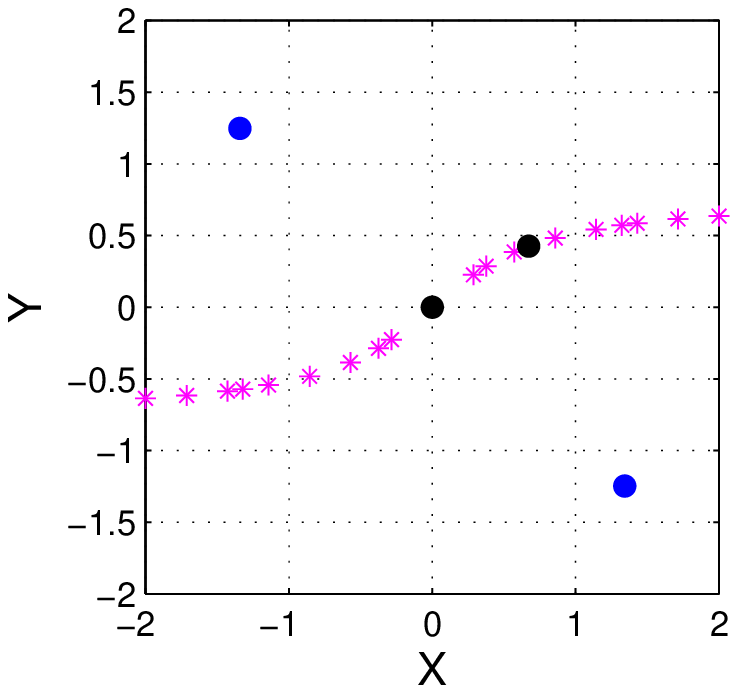} 
		\includegraphics[height=.22\textheight]{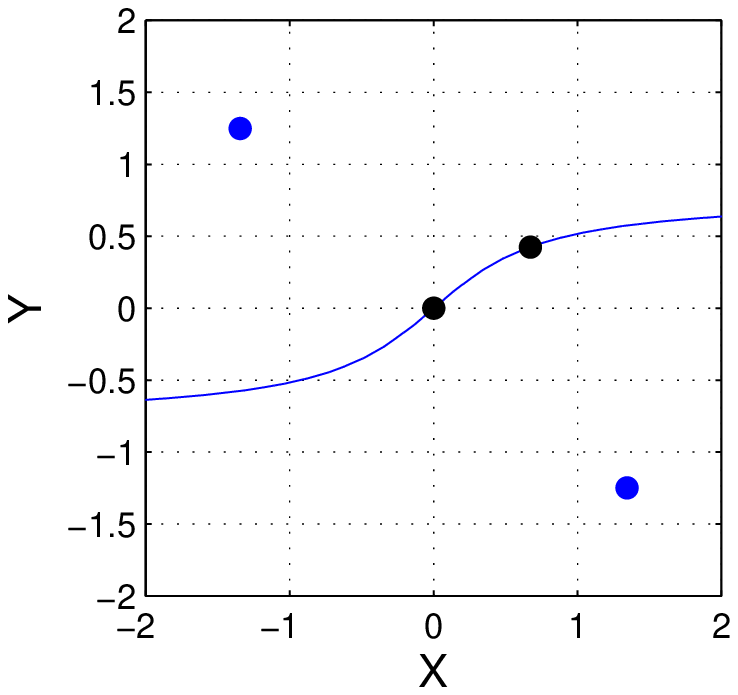} 
		\caption{Set of points lying on  the curve separating the domains of attraction of $E_1$ and $E_2$ (left) and the  reconstruction of the separatrix curve (right). The black and blue circles represent the unstable origin, the coexistence saddle point and the stable equilibria,  respectively.
		}
		\label{curva_ris}
	\end{center}
\end{figure}

To test our detection-interpolation algorithm when bistability occurs we choose the parameters as follows: $r = 0.7895$, $m = 0.7885$, $p = 0.225$, $q = 0.2085$, $Kp = 12$ and $Kq = 10$. With this choice the equilibria $E_1=(10,0)$ and $E_2=(0,16.5)$ are stable, the origin $E_0$ is unstable  and $E_3 \approx (3.8757,3.1919)$ is the saddle coexistence equilibrium point partitioning  the phase space domain of the system \eqref{model2d}. The stable equilibria of \eqref{model2dris} are $E_1^{*}  \approx (1.3436,-1.2482)$, $E_2^{*}  \approx (-1.3436,1.2482)$ and the coexistence saddle point is $E_3^{'} \approx (0.6717,  0.4252)$. In view of the above considerations we can identify $E_1^{*}$ and $E_2^{*}$  with $E_1^{'}=(1,0)$ and $E_2^{'}=(0,1)$.
After translating the problem in the positive plane with the substitutions \eqref{traslaz}, we can apply the detection-interpolation routine.
In this case we choose:  $n=15$, $\gamma=4$, $tol=10^{-4}$, $t=40$,  $\varepsilon=0.1$, $ d^{PU}=3$,   $K=4$.
Figure \ref{curva_ris} shows the separatrix points (left) and   the separatrix curve (right) in the phase plane of the system \eqref{model2dris}.
Using again the transformation \eqref{tras} we obtain the curve separating the basins of attraction of $E_1$, $E_2$, shown in Figure \ref{curve_no_ris} (right).


\end{document}